\documentclass{article}
\font\Bb=msbm10
\def\CC{\hbox{\Bb C}}
\def\NN{\hbox{\Bb N}}
\def\RR{\hbox{\Bb R}}
\def\PP{\hbox{\Bb P}}
\usepackage{amssymb,latexsym,amsmath}

\newtheorem{theorem}{Theorem}[section]

\newtheorem{proposition}[theorem]{Proposition}

\title{Orthogonal polynomials associated with an inverse quadratic spectral transform}
\author{Manuel Alfaro $^{*}$, Ana Pe\~{n}a $^{*}$ , M. Luisa
Rezola \thanks{Partially supported by MEC of Spain under
Grant MTM2009-12740-C03-03 and the DGA project E-64 (Spain).}
\\ Departamento de Matem\'{a}ticas and IUMA \\ Universidad de Zaragoza (Spain)
\and Francisco Marcell\'{a}n \thanks{Partially supported by MEC of Spain under
Grant MTM2009-12740-C03-01 (Spain).} \\ Departamento de Matem\'{a}ticas \\ Universidad Carlos III de Madrid (Spain) }
\date{}

\begin{document}

\maketitle

\bigskip

\begin{abstract}
Let $\{ P_n \}_{n\ge0}$ be a sequence of monic orthogonal polynomials with respect to a quasi--definite linear functional $u$ and  $\{ Q_n \}_{n\ge0}$ a sequence of polynomials defined by
$$Q_n(x)=P_n(x)+s_n\,P_{n-1}(x)+t_n\,P_{n-2}(x),\quad n\ge1,$$
with $t_n \not= 0$ for $n\ge2$.

We obtain a new characterization of the orthogonality of the sequence $\{ Q_n \}_{n\ge0}$ with respect to a linear functional $v$, in terms of the coefficients of a quadratic polynomial $h$ such that $h(x)v= u$.

We also study some cases in which the parameters $s_n$ and $t_n$ can be computed more easily, and give several examples.

Finally, the interpretation of such a perturbation in terms of the Jacobi matrices associated with $\{ P_n \}_{n\ge0}$  and $\{ Q_n \}_{n\ge0}$ is presented.
\end{abstract}

\medskip


\medskip


\bigskip

\noindent {\it Address of the corresponding author}: \medskip
\\Francisco Marcell\'{a}n
\\Departamento de Matem\'{a}ticas \\Universidad Carlos III de Madrid
\\28911 Legan\'{e}s (Spain)\\e.mail: pacomarc@ing.uc3m.es\\
Phone: 34 916249442, \, Fax: 34 916249151.

\newpage

\newpage

\section{Introduction}
\setcounter{equation}{0}

 Let $u$ be a linear functional defined in
the linear space $\mathbb{P}$ of polynomials with complex
coefficients.

The sequence of complex numbers $\{u_n\}_{n\ge0}$ such that $u_n=(u,x^n )$,  is said to be the moment sequence associated with the linear functional $u$. The corresponding $z$-transform $\sum_{n=0}^{\infty} u_n \, z^{-(n+1)}$ defines a function $S(z)$ which is analytic in a neighbourhood of infinity. This function is known in the literature as Stieltjes function associated with the linear functional $u$ (see \cite{M0}).

The linear functional $u$ defines a bilinear functional on $\PP \times \PP$ in the following way: $( p,q )_u =(u,pq)$, for every $p, q \in \PP$. If $H=[(x^k,x^l)_u]_{k,l=0}^{\infty}$ denotes the Gram matrix of the above bilinear functional with respect to the canonical basis of $\PP$, then $H$ is a Hankel matrix, i.e. the entries in every antidiagonal are equal.

The linear functional $u$ is said to be quasi--definite if the leading principal
submatrices $H_n, \, n\ge0 \,,$ of $H$ are non singular. In such a situation, there exists a sequence of monic polynomials
$\{P_n\}_{n\ge0}$ with $\deg P_n=n$ such that $( u,P_nP_m )=k_n \,\delta_{n,m}$ with $k_n\not=0.$
The sequence $\{P_n\}_{n\ge0}$ is said to be the  sequence of monic
orthogonal polynomials (SMOP) with respect to the linear
functional $u.$

Taking into account that the multiplication operator by $x$ is a symmetric operator with respect to the bilinear functional $(.,.)_u$, the SMOP $\{P_n\}_{n\ge0}$ satisfies for $n\ge0$, a three--term recurrence relation $$xP_n(x)=P_{n+1}(x)+\beta_nP_n(x)+\gamma_nP_{n-1}(x),$$ \noindent with $P_0(x)=1$ and $P_{-1}(x)=0$ where
$\gamma_n\not=0$ for every $n \ge1.$ Conversely, if a sequence of monic polynomials $\{P_n\}_{n\ge0}$ satisfies a recurrence
relation like the previous one, then there exists a unique quasi--definite linear
functional $u$ such that  $\{P_n\}_{n\ge0}$ is the corresponding
SMOP. This result is known in the literature as Favard's theorem (see \cite{Ch}).

For the SMOP $\{P_n\}_{n\ge0}$ there exists a sequence $\{\nu_n\}_{n\ge0}$ of linear functionals such that $ (\nu_m, P_n )=\delta_{nm}$. This sequence
$\{\nu_n \}_{n\ge0}$ is said to be the dual sequence of $\{P_n\}_{n\ge0}$. Notice that $\displaystyle \nu_n=\frac{P_n\,u}{(u, P_n^2 )}$, where for every $p \in \PP$, $p\,u$ denotes the linear functional such that $(p\, u, q )= (u, p q )$, $q \in \PP$. The linear functional
$\displaystyle \nu_0=\frac{u}{(u, P_0^2)}$ is said to be the normalization of $u$ in the sense that $ (\nu_0,1 )=1$. Moreover, the linear functional $(x-a)^{-1}u$ is defined on $\PP$ by
\begin{eqnarray*}
((x-a)^{-1} u, P )=(u, \frac{P(x)-P(a)}{x-a}) ,\, P \in \PP.
\end{eqnarray*}

If ${\cal{P}}=(P_0, P_1, \dots )^T$ denotes the column vector associated with an SMOP, then the three--term recurrence relation reads $x{\cal{P}}=J_{\cal{P}}\,{\cal{P}}$ where $J_{\cal{P}}$ is a tridiagonal matrix such that the diagonal entries are $\{\beta_n\}_{n\ge0}$, the upper diagonal entries are $1$, and the subdiagonal entries are $\{\gamma_n\}_{n\ge1}$. It is called a monic Jacobi matrix (see \cite{Ch}).

Thus, for a quasi-definite linear functional $u$ we have three characteristic elements

\begin {itemize}
\item [(i)]  The Stieltjes function.
\item [(ii)] The SMOP $\{P_n\}_{n\ge0}$.
\item [(iii)] The monic Jacobi matrix $J_{\cal{P}}$.
\end{itemize}

The connections between these three basic elements in the theory of orthogonal polynomials with Pad\'e rational approximants, quadrature formula, and spectral theory of symmetric operators have been very fruitful.

In a more general framework, in \cite{X1} and \cite{X2} the author deals with quasi-orthogonal polynomial sequences $\{Q_n\}_{n\ge0}$ of order $r$ with respect to a positive Borel measure $\mu$. They are defined by $Q_{n}(x) = P_{n}(x)+ b_{n,n-1} P_{n-1}(x)+ ...+ b_{n,n-r} P_{n-r}(x)$ with $b_{n,n-r}\neq 0,$ where $\{P_n\}_{n\ge0}$ is the sequence of orthonormal polynomials with respect to $\mu$. When the polynomials $Q_n$ are characteristic polynomials of some symmetric tridiagonal matrices with positive subdiagonal entries, many properties of their zeros can be deduced from such a formulation. Furthermore, a positive quadrature formula with $n$ nodes which is exact for polynomials of degree at most $2n-r-1$, $0\leq r \leq n,$ is based on the zeros of the quasi-orthogonal polynomial $Q_n$ and, as a consequence, for a fixed $n,$ every positive quadrature formula is a Gaussian quadrature formula for some nonnegative measure.

Many authors have been interested in the study of SMOP's associated with some perturbations of linear functionals. In particular, the so called direct problems deal with some canonical perturbations of linear functionals. For instance, the analysis of the parameters of the three-term recurrence relation when a positive definite linear functional is multiplied by a polynomial (Christoffel transformation)  has been done from a numerical point of view in \cite{G}, from the point of view of  the corresponding monic Jacobi matrices in \cite{BuM}, and taking into account the associated Stieltjes functions, in \cite{M} and \cite{SVZ}. The relation between the corresponding SMOP in the case of quasi-definite linear functionals has been deduced in \cite{M3}.
The addition of Dirac masses and their derivatives to a quasi-definite linear functional is called Uvarov transformation (see \cite{GHH}, \cite{M2}, \cite{SVZ}). They play an important role in the spectral analysis of higher order linear differential operators whose eigenfunctions are SMOP (see \cite{GH}). They are related to bispectral problems and self-similar reductions of differential  operators (\cite{M}). The connection between the corresponding monic Jacobi matrices has been studied in \cite{BuM}.

On the other hand, given a quasi-definite linear functional $u$ and polynomials $p,q$ the analysis of the quasi-definiteness of a linear functional $v$ such that $pv=qu$ has been done in the framework of the so-called generalized Christoffel formula when $u$ is a positive definite linear functional and $q/p$ is a positive rational function in the support of $u$ (see \cite{G}). In such a case the connection between the monic Jacobi matrices has been done in \cite{EK}. For the quasi-definite case, the relation between the corresponding SMOP has been given in \cite{M3}.

In the general situation, if $\deg p=1$ and $q$ is a constant, necessary and sufficient conditions for the quasi-definite character of $v$ have been done in \cite{M0}. There, the expression of the perturbed SMOP in term of the initial one as well as the relation between the parameters of the three-term recurrence relation is stated. This perturbation is called Geronimus transformation. The relation between the corresponding monic Jacobi matrices is analyzed in \cite{BuM} and \cite{M2}. Finally, for the Stieltjes functions associated with $u$ and $v$ the connection is given in \cite{SVZ}.

When $\deg p=\deg q=1$, the relation between the SMOP associated with $u$ and $v$ is stated in \cite{AMPR2}.

Eventually, when $p(x)=x^2$ and $q$ is a constant, in \cite{BeM} the authors focus the attention in the invariance of the semiclassical character for such a transformation.

All these cases correspond to the so called direct problems, i.e., how a perturbation in a linear functional concerns to Stieltjes functions, monic Jacobi matrices, and orthogonal polynomials, respectively.

The study of inverse problems, i.e. to find the relation between two linear functionals assuming their corresponding SMOP's are related in terms of an algebraic expression, that it is said to be of type k-j, $Q_{n}(x)+ a_{n,n-1} Q_{n-1}(x)+ ...+ a_{n,n-k+1} Q_{n-k+1}(x) = P_{n}(x)+ b_{n,n-1} P_{n-1}(x)+ ...+ b_{n,n-j+1} P_{n-j+1}(x)$ has been firstly considered by Geronimus in \cite{Ge}  in the case 1-2. Later on, the relation 1-3 has been  analyzed in  \cite{BM} as well as in \cite{HHR} from the point of view of the corresponding parameters of the three--term recurrence relation. In this last contribution the authors also deduce the second order linear differential equation satisfied by the perturbed sequence of orthogonal polynomials when the initial sequence is a classical one. We also point out that in \cite{WK1} and  \cite{WK2} the authors focus the attention in the study of these perturbed families when the initial orthogonal polynomial sequences are the Chebyshev polynomials of first and second kind, respectively. In particular, the integral representation of the corresponding linear functional is obtained. Such a kind of relations appear in the framework of coherent pairs of linear functionals related to Sobolev inner products. This problem has been pointed out in \cite{MP}. Relations of type 2-2 have been completely studied in \cite{AMPR1}. In a very recent paper \cite{AMPR3}, the authors solved a problem studied in a particular case by Grinshpun (see \cite{Gr}), concerning the orthogonality of sequences of monic polynomials $\{Q_n\}_{n\ge0}$ with $Q_n(x)=P_{n}(x)+s P_{n-1}(x)+t P_{n-2}(x)$ where $\{P_n\}_{n\ge0}$ is a given SMOP.

In \cite{AMPR1} and \cite{MP}, the orthogonality of the sequence of monic polynomials $Q_n$  is given in terms of the some constant sequences. Moreover, these constants are those appearing in the polynomial which relates the corresponding linear functional associated with  $\{P_n\}_{n\ge0}$ and $\{Q_n\}_{n\ge0}$. More precisely, in \cite{MP} for a case 1-2 the authors proved the following results:

Let $\{ P_n \}_{n\ge0}$ be an SMOP with respect to a quasi-definite linear functional $u$ with recurrence coefficients $\{\beta_n\}_{n\ge0}$ and $\{\gamma_n\}_{n\ge1}$. We define a sequence of monic polynomials  $\{ Q_n \}_{n\ge0}$ by
$Q_n = P_n+ \mu_n P_{n-1}, \, n\ge1$, where $\mu_n$ are complex numbers and $\mu_1 \not= 0$. The sequence $\{ Q_n \}_{n\ge0}$ is an SMOP with respect to a quasi-definite linear functional $v$ if and only
$P_n^* (x_1; x_1-\beta_0+\mu_1,0)\not=0$ for all $n\ge1$ and $\mu_n= - \frac{P_n^* (x_1; x_1-\beta_0+\mu_1,0)}{P_{n-1}^* (x_1; x_1-\beta_0+\mu_1,0)}$, where $\{P_n^* (x; x_1-\beta_0+\mu_1,0)\}_{n \ge 0}$ denotes the co-recursive sequence at level zero corresponding to the sequence $\{ P_n \}_{n \ge 0}$ (see \cite {Ch}, \cite {MDR}), and $x_1 =\beta_1-\mu_2-\frac{\gamma_1}{\mu_1}$.

This condition is equivalent to the fact that for $n\ge1$ the parameters $\mu_n$ do not vanish and satisfy
\begin{equation} \label{Petro}
\beta_n - \mu_{n+1}-\frac{\gamma_n}{\mu_n}= x_1, \quad n \ge 1 .
\end{equation}
Moreover, $(x-x_1)v = (\beta_0-x_1-\mu_1) u$. (See \cite[Theorem 2]{MP}). This result will be used many times in this paper.

In the same way in \cite{AMPR1}, the authors studied the relations 2-2 ($P_n+s_nP_{n-1}=Q_n+t_nQ_{n-1}$). They obtained the characterization of the orthogonality of the polynomials $Q_n$ in terms of the existence of two constant sequences. Besides, these constants $a,\tilde{a}$ are the values  appearing in the relation between the two linear functionals, i.e. $(x-\tilde{a})u=k(x-a)v$,  ($a,\tilde{a},k \in \CC$).

A natural question is: Do these constant sequences depend on the coefficients of the polynomial relating the two linear functionals for more general cases?

The aim of our  contribution is twofold. First,  we find necessary and sufficient conditions that the sequence $\{Q_n\}_{n\ge0}$, with $Q_n(x)=P_{n}(x)+s_n P_{n-1}(x)+t_n P_{n-2}(x)$, must satisfy in order to be orthogonal with respect to a linear functional $v$. Indeed,  if $u$ is the linear functional such that $\{ P_n \}_{n\ge0}$ is the corresponding SMOP, then there exists a monic quadratic polynomial $h(x)= x^2+ax+b$ such that $h(x)v=k \, u$. We prove that in this case, there exist two constant sequences and the values of these constants $a$ and $b$ are precisely the coefficients and not the zeros of the quadratic polynomial  $h(x)=x^2+ax+b$ such that $h(x)v= k \, u$. A connection between  the corresponding monic Jacobi matrices is stated.

Second, we obtain necessary and sufficient conditions in order to the above relation 1-3 can be decomposed in two relations 1-2 and then proceed by iteration. This is related with a two step Geronimus transformation. In this case, we prove that the constant sequences can be decomposed in another two simplest constant sequences, associated with the zeros of the polynomial $h(x)$. Also, as a consequence, a matrix interpretation using $LU$  and $UL$ factorization is done. Finally, when $\{P_n\}_{n\ge0}$ is a symmetric  SMOP, i.e. $ P_n(-x)=(-1)^n P_n(x)$ and we assume $s_n=0$, then we prove that the 1-3 relation yields two 1-2 relations.

The structure of the manuscript is the following. In Section 2 we present the basic background and we give a new characterization of the orthogonality of the sequence $\{Q_n\}_{n\ge0}$ (Theorem  \ref {mainresult}) in terms of some constant relations for $\{\beta_n\}, \, \{\gamma_n\}, \, \{s_n\}$, and $ \{t_n\}$. They have not yet studied in the literature. In Section 3 we deduce necessary and sufficient conditions for the split of a 1-3 relation in two 1-2 relations. The matrix interpretation of this problem in terms of monic Jacobi matrices is done. Some illustrative examples are considered.

\section{A new characterization of orthogonality}
\setcounter{equation}{0}

From now on, $\{ P_n \}_{n \ge 0}$ will denote an SMOP with respect to a quasi-definite linear functional $u$ which satisfies
the three--term recurrence relation
\begin{align} \label{recurrenciaPn}
P_{n+1}(x)&=(x-\beta_n)P_n(x)-\gamma_nP_{n-1}(x), \quad n\ge0,
\\P_0(x)&=1, \quad P_{-1}(x)=0,\nonumber
\end{align}
\noindent{where} $\{\beta_n\}_{n\ge0}$ and $\{\gamma_n\}_{n\ge1}$
are sequences of complex numbers with $\gamma_n \not=0$ for $n \ge
1$.

Given two sequences of complex numbers $\{s_n\}_{n\ge1}$ and $\{t_n\}_{n\ge2}$  we define for $n\ge1$ the sequence of monic polynomials $ \{ Q_n \}_{n \ge 0}$ such that
\begin{equation}\label{relacionuno-tres}
Q_n(x)=P_n(x)+s_nP_{n-1}(x)+t_nP_{n-2}(x), \quad  \text {with} \quad t_n \not= 0, \, n\ge2.
\end{equation}

\medskip

First, we characterize the orthogonality of the sequence $\{Q_n \}_{n \ge 0}$.

\begin{proposition}\label{caracterizacion ortogonalidad} Let $\{ P_n \}_{n \ge 0}$ an SMOP with recurrence
coefficients $\{\beta_n\}_{n\ge0}$ and $\{\gamma_n\}_{n\ge1}$. We define a sequence
$\{Q_n\}_{n \ge 0}$ of monic polynomials by formula
(\ref{relacionuno-tres}), i.e.,
\begin{equation*}
Q_n(x)=P_n(x)+s_nP_{n-1}(x)+t_nP_{n-2}(x),\quad n\ge1,
\end{equation*}
\noindent where  $s_n$ and $t_n$ are complex numbers with $t_n\not=0$,
for all $n \ge 2$. Then
$\{ Q_n \}_{n\ge0}$ is an SMOP with recurrence coefficients
$\{\tilde{\beta}_n\}_{n\ge0}$ and $\{\tilde{\gamma}_n\}_{n\ge1}$,  $\tilde{\gamma}_n \not= 0$, if and only if the following
formulas hold:
\begin{equation}\label{uno}
s_{n-1}\tilde{\gamma}_n=s_n\gamma_{n-1}+t_n(\beta_{n-2}-\tilde{\beta}_n),
\quad n \ge 2 \,,
\end{equation}
\begin{equation}\label{dos}
t_{n-1}\tilde{\gamma}_n= t_n\gamma_{n-2}, \quad n \ge 3 \,,
\end{equation}
\noindent where the coefficients $\tilde{\beta}_n$ and
$\tilde{\gamma}_n$ are defined by
\begin{equation}\label{betatilden}
\tilde{\beta}_n=\beta_n+s_n-s_{n+1}, \quad n \ge 0 \,,
\end{equation}
\begin{equation}\label{gammatilden}
\tilde{\gamma}_n=\gamma_n+t_n-t_{n+1}+s_n(\beta_{n-1}-\tilde{\beta}_n),
\quad n \ge 1,
\end{equation}
\noindent with $s_0=t_0=t_1=0.$
\end{proposition}

\textbf{Proof.} Inserting formula (\ref{recurrenciaPn}) in
(\ref{relacionuno-tres}) and applying (\ref{relacionuno-tres}) to
$xP_n(x)$, we get

\begin{align*}
Q_{n+1}(x)&=xQ_n(x)+(s_{n+1}-s_n-\beta_n)P_n(x) \\
&+(t_{n+1}-\gamma_n-s_n\beta_{n-1}-t_n)P_{n-1}(x)\\
&-(s_n\gamma_{n-1}+t_n\beta_{n-2})P_{n-2}(x)-t_n\gamma_{n-2}P_{n-3}(x),
\quad n \ge 1\,,
\end{align*}
provided we substitute there $xP_{n-1}(x)$ and
$xP_{n-2}(x)$, using again (\ref{recurrenciaPn}). Now, formula
(\ref{relacionuno-tres}) applied to $P_n(x)$ and the definition of
$\tilde{\beta}_n$ (see (\ref{betatilden})), yield

\begin{align*}
Q_{n+1}(x)&=(x-\tilde{\beta}_n) Q_n(x) \\
&+\left[(t_{n+1}-\gamma_n-t_n-s_n(s_{n+1}-s_n-\beta_n+\beta_{n-1})\right]P_{n-1}(x)\\
&-\left[s_n\gamma_{n-1}+t_n(s_{n+1}-s_n-\beta_n+\beta_{n-2})\right]P_{n-2}(x)\\
&-t_n\gamma_{n-2}P_{n-3}(x), \quad n \ge 1 \,.
\end{align*}
So $\{ Q_n \}_{n\ge0}$ is an SMOP if and only if
there exists a sequence of complex numbers
$(\tilde{\gamma}_n)_{n\ge1}$, with $\tilde{\gamma}_n \not=0$ for
$n\ge1$, such that
\begin{align*}
&\left[ t_{n+1}-\gamma_n-t_n-s_n(s_{n+1}-s_n-\beta_n+\beta_{n-1}) \right] P_{n-1}(x)
\\
&-\left[s_n \gamma_{n-1} + t_n(s_{n+1}-s_n-\beta_n + \beta_{n-2})
\right]  P_{n-2}(x) - t_n \gamma_{n-2} P_{n-3}(x) \\
&=- \tilde{\gamma}_nQ_{n-1}(x).
\end{align*}

\noindent{Furthermore}, $\tilde{\beta}_n$ and $\tilde{\gamma}_n$ are
the coefficients of the three--term recurrence relation for $Q_n.$

Using again (\ref{relacionuno-tres}), straightforward calculations yield
(\ref{uno}) and (\ref{dos}).$\quad \Box$

\medskip
This result can be also seen in \cite{BM}.

\medskip

\noindent \textbf {Remark.}  Observe that the conditions $\tilde{\gamma}_n \not= 0, \, \forall n\ge 3 \,,$ follow from (\ref{dos}).

\bigskip

If the sequence $\{ Q_n \}_{n \ge 0}$ given by (\ref{relacionuno-tres}) is an SMOP with respect to a quasi-definite linear functional $v$, then it is well known that  the linear functional $u$ is a quadratic polynomial modification of $v$, that is, $h(x) v = k \,  u$ with $h(x) = x^2 +ax+b$ and $k$ a nonzero complex number.

In the next theorem, we  show that the orthogonality of the sequence $\{Q_n \}_{n \ge 0}$ can be characterized by the fact that there are two sequences depending on the parameters $s_n, t_n, \beta_n, \gamma_n$ which remain constant and these constants are precisely the coefficients $a$ and $b$ of the polynomial $h$.

\bigskip

\begin{theorem}\label{mainresult}
Let $\{ P_n \}_{n \ge 0}$ be
an SMOP with recurrence
coefficients $\{\beta_n\}_{n\ge0}$ and $\{\gamma_n\}_{n\ge1}$. We define a sequence
$\{Q_n\}_{n \ge 0}$ of monic polynomials by the formula
(\ref{relacionuno-tres}), i.e.,
\begin{equation*}
Q_n(x)=P_n(x)+s_nP_{n-1}(x)+t_nP_{n-2}(x),\quad n\ge1,
\end{equation*}
\noindent where  $s_n$ and $t_n$ are complex numbers with $t_n\not=0$,
for all $n \ge 2$ and $s_0 = t_0 = t_1 = 0$.
Then the following statements are equivalent:
\item [(i)] $\{ Q_n \}_{n\ge0}$ is an SMOP.
\item [(ii)] There exist two constant sequences $\{ A_n \}_{n \ge 1}$ and $\{ B_n \}_{n \ge 1}$ such that, for $n \ge 1$,
\begin{align}\label{a}
& A_n = \frac{s_n}{t_{n+1}} [\gamma_{n+1} + t_{n+1} - t_{n+2} + s_{n+1}(\beta_n - {\beta}_{n+1} - s_{n+1}
+ s_{n+2})]  \\
&+ s_{n+1} - \beta_{n-1} - {\beta}_n = a ,\nonumber
\end{align}
\begin{align}\label{b}
& B_n = \frac{1}{t_{n+1}}[\gamma_{n+1} - t_{n+2} + s_{n+1}(\beta_n - \beta_{n+1} - s_{n+1} +
s_{n+2})] \\
&\times [\gamma_n + {t_n} - t_{n+1} + s_n (s_{n+1}-\beta_n)] + {t_n} -
\gamma_{n-1}+(s_{n+1}-\beta_n)(s_n-\beta_{n-1}) = b,\nonumber
\end{align}
where, by convention, we take $\gamma_0=0$. Besides, the values of $\tilde{\gamma}_1$ and $\tilde{\gamma}_2$ in (\ref{gammatilden}) have to be different from zero.

Furthermore, if $u$ and $v$ are the linear functionals associated with the sequences
$\{ P_n \}_{n\ge0}$ and $\{ Q_n \}_{n\ge0}$, respectively, normalized by
$\langle u,1 \rangle = \langle v,1 \rangle = 1,$ then
\begin{equation} \label{relacionfun}
 (x^2 +a x + b) v = k \,  u .
\end{equation}
\noindent with $k \in \CC \setminus \{0 \}$ .
\end{theorem}

\textbf{Proof.}
For simplicity we write the proof of this theorem in terms of $\tilde{\beta}_n$ and $\tilde{\gamma}_n$ and then it is enough to use (\ref{betatilden}) and (\ref{gammatilden}).

$(i)\Rightarrow (ii):$

If $\{Q_n\}_{n \ge 0}$ is an SMOP, by Proposition
\ref{caracterizacion ortogonalidad}, (\ref{uno}) and (\ref{dos}) hold. Inserting in (\ref{uno}) the expression of $\gamma_{n-1}$ in terms of $\tilde{\gamma}_{n+1}$ given by (\ref{dos}) we get for $n \ge 2$,
\begin{equation}\label{intermedia}
\frac{s_{n-1}}{t_n} \tilde{\gamma}_n -\tilde{\beta}_{n-2}-\tilde{\beta}_{n-1} +s_{n-2}
 = \frac{s_n}{t_{n+1}} \tilde{\gamma}_{n+1} -\tilde{\beta}_{n-1}-\tilde{\beta}_n +s_{n-1}
\end{equation}
and then we get (\ref{a}).

Now, we will deduce (\ref{b}).

From definition of $\tilde{\gamma}_n$, using (\ref{dos}) we have
\begin{equation*}
\tilde{\gamma}_n = \frac{t_{n+1}}{t_{n+2}} \tilde{\gamma}_{n+2} + t_n - t_{n+1} +
s_n (\tilde{\beta}_{n-1} - s_{n-1} + s_n -\tilde{\beta}_n) \,, \quad n \ge 1\,.
\end{equation*}
Multiplying in the above expression by $\tilde{\gamma}_{n+1} / t_{n+1}$
\begin{equation*}
\frac{\tilde{\gamma}_n\tilde{\gamma}_{n+1}}{t_{n+1}}+(s_{n-1} - \tilde{\beta}_{n-1})\frac{s_n \tilde{\gamma}_{n+1}}{t_{n+1}} =
\frac{\tilde{\gamma}_{n+1} \tilde{\gamma}_{n+2}}{t_{n+2}} +
(s_n - \tilde{\beta}_n) \frac{s_n \tilde{\gamma}_{n+1}}{t_{n+1}} + \left( \frac{t_n}{t_{n+1}} - 1 \right) \tilde{\gamma}_{n+1} \,.
\end{equation*}
Using (\ref{intermedia}), for $n + 1$ instead of $n$, in the expression $\displaystyle \frac{s_n \tilde{\gamma}_{n+1}}{t_{n+1}}$ which appears in the right hand side of the above formula, for $n \ge 1$ we obtain
\begin{align*}
&\frac{\tilde{\gamma}_n\tilde{\gamma}_{n+1}}{t_{n+1}} + (s_{n-1} - \tilde{\beta}_{n-1}) \left( \frac{s_n \tilde{\gamma}_{n+1}}{t_{n+1}} - \tilde{\beta}_n \right) \\
&= \frac{\tilde{\gamma}_{n+1}\tilde{\gamma}_{n+2}}{t_{n+1}} + (s_n - \tilde{\beta}_n ) \left( \frac{s_{n+1} \tilde{\gamma}_{n+2}}{t_{n+2}} - \tilde{\beta}_{n+1} \right) \\
&+ \left(\frac{t_n}{t_{n+1}} - 1 \right) \tilde{\gamma}_{n+1} + s_n (\tilde{\beta}_{n-1} + s_n - s_{n-1} - \tilde{\beta}_n)  \,.
\end{align*}

Besides, from (\ref{dos}) and according to the definition of $\tilde{\gamma}_{n-1}$ we have, for $n \ge 1$,
\begin{equation*}
\frac{t_n}{t_{n+1}} \tilde{\gamma}_{n+1} = \tilde{\gamma}_{n-1} - t_{n-1} + t_n -
s_{n-1} (\tilde{\beta}_{n-2} + s_{n-1} - s_{n-2} -\tilde{\beta}_{n-1})
\end{equation*}
and, therefore,
\begin{align}\label{intermedia3}
&\frac{\tilde{\gamma}_n\tilde{\gamma}_{n+1}}{t_{n+1}} + (s_{n-1} - \tilde{\beta}_{n-1}) \left(\frac{s_n \tilde{\gamma}_{n+1}}{t_{n+1}} - \tilde{\beta}_n \right)\nonumber \\
&+ t_{n-1} - \tilde{\gamma}_{n-1} - \tilde{\gamma}_n +
s_{n-1} (\tilde{\beta}_{n-2} + s_{n-1} - s_{n-2} -\tilde{\beta}_{n-1}) \\\nonumber
&= \frac{\tilde{\gamma}_{n+1}\tilde{\gamma}_{n+2}}{t_{n+1}} + (s_n - \tilde{\beta}_n) \left( \frac{s_{n+1} \tilde{\gamma}_{n+2}}{t_{n+2}} - \tilde{\beta}_{n+1} \right) \\ \nonumber
&+ t_n - \tilde{\gamma}_n - \tilde{\gamma}_{n+1} +
s_n (\tilde{\beta}_{n-1} + s_n - s_{n-1} -\tilde{\beta}_n). \nonumber
\end{align}

Thus, by straightforward computations, using (\ref{betatilden}) and (\ref{gammatilden}) we obtain (\ref{b}).

$(ii)\Rightarrow (i):$

First, we will prove that (\ref{intermedia}) and (\ref{intermedia3}) yield (\ref{dos}). Taking into account the new expression of $\displaystyle \frac{s_{n+1} \tilde{\gamma}_{n+2}}{t_{n+2}}$, obtained from (\ref{intermedia}) written for $n + 1$ instead of $n$, we can reformulate (\ref{intermedia3})
\begin{align*}
&\frac{\tilde{\gamma}_n\tilde{\gamma}_{n+1}}{t_{n+1}} + (s_{n-1} - \tilde{\beta}_{n-1}) \left( \frac{s_n \tilde{\gamma}_{n+1}}{t_{n+1}}  - \tilde{\beta}_n \right) \\
&+ t_{n-1} - \tilde{\gamma}_{n-1} +
s_{n-1} (\tilde{\beta}_{n-2} + s_{n-1} - s_{n-2} -\tilde{\beta}_{n-1}) \\
&= \frac{\tilde{\gamma}_{n+1}\tilde{\gamma}_{n+2}}{t_{n+1}} + \left(s_n - \tilde{\beta}_n \right) \left[\frac{s_n \tilde{\gamma}_{n+1}}{t_{n+1}} - \tilde{\beta}_{n-1} + s_{n-1} - s_n\right]\\
&+ t_n  - \tilde{\gamma}_{n+1} +
s_n (\tilde{\beta}_{n-1} + s_n - s_{n-1} -\tilde{\beta}_n)  \,.
\end{align*}
Using the definition of $\tilde{\gamma}_{n-1}$ in the left hand side of the above expression and simplifying we get

\begin{equation*}
\frac{\tilde{\gamma}_{n+1}}{t_{n+1}} \left[\tilde{\gamma}_n + s_n(\tilde{\beta}_n - s_n + s_{n-1}-\tilde{\beta}_{n-1}) + t_{n-1} \right] = \gamma_{n-1} + \frac{\tilde{\gamma}_{n+1}\tilde{\gamma}_{n+2}}{t_{n+2}} \,.
\end{equation*}
Using again the definition of $\tilde{\gamma}_n$ we have
\begin{equation*}
\frac{\tilde{\gamma}_{n+1}}{t_{n+1}} (\gamma_n + t_n) = \gamma_{n-1} + \frac{\tilde{\gamma}_{n+1}\tilde{\gamma}_{n+2}}{t_{n+2}} \,,
\end{equation*}
that is, for $n \ge 1$
\begin{equation}\label{gammasintilde}
\frac{\tilde{\gamma}_{n+1}}{t_{n+1}} \left(\gamma_n - \frac{t_{n+1}}{t_{n+2}}\tilde{\gamma}_{n+2}\right)
= \gamma_{n-1} + \frac{t_{n+1}}{t_{n+2}} \tilde{\gamma}_{n+1}\,.
\end{equation}

Since $\gamma_0 - \frac{t_1}{t_2}\tilde{\gamma}_2 = 0$, we deduce $\gamma_n
= \frac{t_{n+1}}{t_{n+2}}\tilde{\gamma}_{n+2}$, for $n \ge 1$, and therefore (\ref{dos}) holds. Finally, we
point out that from (\ref{intermedia}) and (\ref{dos}) we can derive (\ref{uno}).

To conclude the proof it remains to deduce the relation between the functionals $u$ and $v$ in terms of the
constants $a$ and $b$.

If we expand the linear functional $u$ in the dual basis $\displaystyle \left\{ \frac{Q_jv}{\langle v , Q_j^2 \rangle} \right\}_{j \ge 0 }$ of the polynomials $\{Q_j\}_{j\geq0}$ (see \cite{M2}) and taking into account
(\ref{relacionuno-tres}), then
\begin{equation*}
u = \sum_{j=0}^2 \frac{\langle u , Q_j \rangle}{\langle v , Q_j^2 \rangle} Q_jv =
\left( \frac{t_2}{\tilde{\gamma}_1 \, \tilde{\gamma}_2} Q_2 + \frac{s_1}{\tilde{\gamma}_1} Q_1 + 1 \right)
\, v \,.
\end{equation*}
Introducing the explicit expression of the polynomials $Q_2$ and $Q_1$ given by the three--term recurrence
relation we obtain
\begin{equation*}
u =  \frac{t_2}{\tilde{\gamma}_1\tilde{\gamma}_2} \left[(x^2 + (\frac{s_1\tilde{\gamma}_2}{t_2} -
\tilde{\beta}_0 - \tilde{\beta}_1) x + \left(\frac{\tilde{\gamma}_1\tilde{\gamma}_2}{t_2} -
\tilde{\beta}_0\frac{s_1\tilde{\gamma}_2}{t_2} - \tilde{\gamma}_1 + \tilde{\beta}_0\tilde{\beta}_1 \right)
\right] v \,.  \,\, \Box
\end{equation*}

\noindent \textbf {Remark.} The proof of the theorem has been written in terms of $\tilde{\beta}_n$ and $\tilde{\gamma}_n$ not only for simplicity. If we work in a similar way as above but trying to look for that all the
expressions appearing therein depend only on $\beta_n$ and $\gamma_n$, then we can obtain similar
formulas for $a$ and $b$ but a problem appears when we want to come back. Indeed, instead of
(\ref{gammasintilde}) we obtain
\begin{equation*}
\frac{\gamma_{n-1}}{t_n}
\left[\tilde{\gamma}_n - t_n \frac{\tilde{\gamma}_{n-2}}{t_{n-1}}\right]  = \tilde{\gamma}_{n+1} +
t_{n+1}\frac{\gamma_{n-1}}{t_n},\quad n\ge3 \,,
\end{equation*}
and so we can not achieve (\ref{dos}). The main reason is that the trivial extension of (\ref{dos}) to $n = 2 \, (t_1 = \gamma_0 = 0)$ allows us to obtain the value of $\gamma_0$ in terms of $\tilde{\gamma}_2$ as $\gamma_0 =  t_1 \frac{\tilde{\gamma}_2}{t_2}$ but not conversely, i.e. $\tilde{\gamma}_2$ in terms of $\gamma_0$.

\medskip

In the sequel, we present a matrix interpretation of these results in terms of the monic Jacobi matrices associated with the SMOP's $\{P_n\}_{n \ge 0}$  and  $\{Q_n\}_{n \ge 0}$, respectively.

Let ${\cal{P}}=(P_0, P_1, \dots )^T$ and ${\cal{Q}}=(Q_0, Q_1, \dots )^T$ be the column vectors associated with these orthogonal families, and $J_{\cal{P}}$ and $J_{\cal{Q}}$ the corresponding monic Jacobi matrices. Then, the three--term recurrence relations for such SMOP's read $x{\cal{P}}=J_{\cal{P}}\,{\cal{P}}$ and $x{\cal{Q}}=J_{\cal{Q}}\,\cal{Q}$.

Next, we will describe a method to find the matrix $J_{\cal{Q}}$ using the matrix $J_{\cal{P}}$ and the polynomial $h(x)=x^2 +a x + b$.

Taking into account ${\cal{Q}}={\cal{M}}\cal{P}$ where ${\cal{M}}=(m_{k,j})$ is a banded lower triangular matrix such that
$m_{k,k}=1$ and $m_{k,j}=0$ for $k-j>2$, then
$x{\cal{M}}{\cal{P}}=J_{\cal{Q}}{\cal{M}}\cal{P}$
 and, as a consequence, $J_{\cal{P}}{\cal{P}}={\cal{M}}^{-1}J_{\cal{Q}}{\cal{M}}\cal{P}$. Thus, we get
\begin{equation*}
{\cal{M}}J_{\cal{P}}=J_{\cal{Q}}{\cal{M}}.
\end{equation*}

On the other hand, from the classical Christoffel formula (see \cite{G}) we can express $h(x){\cal{P}}$ using the matrix representation
\begin{equation*}
h(x){\cal{P}}={\cal{N}}\cal{Q}
\end{equation*}
\noindent where ${\cal{N}}$ is a banded upper triangular matrix such that $n_{k,k+2}=1$ and $n_{k,j}=0$ for $j-k>2$.
Thus $h(x){\cal{P}}={\cal{N}}{\cal{M}}\cal{P}$, and then
\begin{equation}\label{JP2JPINM}
J_{\cal{P}}^2+aJ_{\cal{P}}+bI={\cal{N}}{\cal{M}}.
\end{equation}

But, from $h(x){\cal{Q}}={\cal{M}}{\cal{N}}\cal{Q}$,  we get
\begin{equation}\label{JQ2JQIMN}
J_{\cal{Q}}^2+aJ_{\cal{Q}}+bI={\cal{M}}{\cal{N}}.
\end{equation}
As a conclusion, we can summarize our process as follows

\medskip

\noindent Step 1: Given $J_{\cal{P}}$, we find the polynomial matrix $J_{\cal{P}}^2+aJ_{\cal{P}}+bI$.

\noindent Step 2: From ${\cal{M}}$ and (\ref{JP2JPINM}) we find ${\cal{N}}$.

\noindent Step 3: From (\ref{JQ2JQIMN}) we obtain the polynomial matrix $J_{\cal{Q}}^2+aJ_{\cal{Q}}+bI$.

\noindent Step 4: Taking into account $J_{\cal{Q}}$ is a tridiagonal matrix, from step 3 we can deduce $J_{\cal{Q}}$, since $\left(J_{\cal{Q}} +\frac{a}{2}\,I\right)^2 = {\cal{M}}{\cal{N}} - \left(b - \frac{a^2}{4}\right)\,I$.

\medskip

A similar analysis for the truncated matrices yields the principal submatrices of $J_{\cal{Q}}$ can be deduced from a rank--one perturbation of the corresponding leading principal submatrices of $J_{\cal{P}}$.

Let denote $({\cal{P}})_n=(P_0, P_1, \dots , P_n)^T$ and $({\cal{Q}})_n=(Q_0, Q_1, \dots , Q_n)^T$. Then the corresponding three--term recurrence relations read
\begin{equation}\label{paco1}
x\, ({\cal{P}})_n = (J_{\cal{P}})_{n+1} ({\cal{P}})_n + P_{n+1} e_{n+1} \,,
\end{equation}
\begin{equation}\label{paco2}
x\, ({\cal{Q}})_n = (J_{\cal{Q}})_{n+1} ({\cal{Q}})_n + Q_{n+1} e_{n+1} \,,
\end{equation}
where the symbol $(A)_n$ stands for the truncation of any infinite matrix $A$ at level $n$, that is, $(J_{{\cal{P}}})_{n+1}$ denotes the leading principal submatrrix of $J_{\cal{P}}$ of size $(n+1) \times (n+1)$, and  $e_{n+1} = (0,\dots ,0 , 1)^T \in \RR^{n+1}$.

Taking into account
\begin{equation}\label{paco3}
({\cal{Q}})_n = ({\cal{M}})_{n+1} ({\cal{P}})_n  \,,
\end{equation}
\noindent and replacing (\ref{paco3}) in (\ref{paco2}) we get
\begin{align*}
&x \, ({\cal{M}})_{n+1} ({\cal{P}})_n = (J_{\cal{Q}})_{n+1} ({\cal{M}})_{n+1} ({\cal{P}})_n + (P_{n+1} + s_{n+1} P_n + t_{n+1} P_{n-1}) e_{n+1}  \\
&= [(J_{\cal{Q}})_{n+1} ({\cal{M}})_{n+1}  + e_{n+1}( s_{n+1} e_{n+1}^T + t_{n+1}d_{n+1}^T)]({\cal{P}})_n + P_{n+1} e_{n+1} \,,
\end{align*}
where $d_{n+1} = (0,\dots ,0 , 1, 0)^T \in \RR^{n+1}$.

Thus, from (\ref{paco1})
\begin{equation*}
({\cal{M}})_{n+1} (J_{\cal{P}})_{n+1} =  (J_{\cal{Q}})_{n+1} ({\cal{M}})_{n+1}  + e_{n+1}(s_{n+1} e_{n+1}^T + t_{n+1} d_{n+1}^T) \,.
\end{equation*}

So, in order to obtain $(J_{\cal{Q}})_{n+1}$ we proceed in three steps.

\medskip

\noindent Step 1: Obtain $({\cal{M}})_{n+1} (J_{\cal{P}})_{n+1}$.

\noindent Step 2: Subtract to the previous one the rank--one matrix $e_{n+1}(s_{n+1} e_{n+1}^T + t_{n+1}d_{n+1}^T)$. This leads to a matrix $G_{n+1}$.

\noindent Step 3: $(J_{\cal{Q}})_{n+1}   = G_{n+1} [({\cal{M}})_{n+1}]^{-1}$.

In the case of quasi-orthogonal polynomials associated with positive Borel measures, this approach has been considered in \cite{X1} and \cite{X2} for symmetric tridiagonal matrices. Nevertheless, not every quasi-orthogonal polynomial can be represented as the characteristic polynomial of a symmetric tridiagonal matrix. Indeed, in \cite[Theorem 3.3]{X2}  a sufficient condition for such a representation is given.

\section {Reducible cases}
\setcounter{equation}{0}

Let  $\{P_n \}_{n \ge 0}$ be an SMOP with respect to a quasi-definite linear functional $u$, and $\{Q_n \}_{n \ge 0}$ a sequence of monic polynomials strictly quasi-orthogonal of order two, satisfying the relation (\ref{relacionuno-tres}).

As a consequence of the Proposition \ref{caracterizacion ortogonalidad} we have that $\{Q_n \}_{n \ge 0}$ is an SMOP if and only if the parameters $s_n$ and $t_n$ satisfy, for $n \ge 3$,

 \begin{equation*}
t_{n+1} = t_n (1- \frac{\gamma_{n-2}}{t_{n-1}})+ \gamma_n +s_n(\beta_{n-1}-\beta_n+ s_{n+1}-s_n) \,,
\end{equation*}
\begin{equation*}
s_{n+1} = -\frac{s_n}{t_n} \gamma_{n-1}+\beta_n-\beta_{n-1}+s_3+\frac{s_2}{t_2}\gamma_1 -\beta_2- \beta_1 \,,
\end{equation*}
\noindent with some initial conditions.

In this section we will study two particular situations where the computations of these parameters can be simplified.

\subsection{Iterative case}

Notice that if $\{Q_n \}_{n\ge0}$ satisfies (\ref{relacionuno-tres}), then the polynomials $Q_n$ cannot be represented as a linear combination of the two consecutive polynomials $P_n$ and $P_{n-1}$. A natural question arises: Can the SMOP $\{ Q_n \}_{n\ge0}$  be generated from $\{ P_n \}_{n\ge0}$ in two steps with the help of an intermediate SMOP? In other words: Do exist two sequences of nonzero complex numbers $\{ \lambda_n \}_{n\ge1}$ and $\{ \mu_n \}_{n\ge1}$ such that $\{ R_n \}_{n\ge0}$ defined by $R_n = P_n+ \mu_n P_{n-1}$ is an SMOP and $Q_n = R_n+ \lambda_n R_{n-1}$? The interest of this question is to simplify the computation of the parameters $s_n$ and $t_n$, in this case via the parameters $\lambda_n$ and $\mu_n$.

If the answer is positive, then we will say that the representation
(\ref{relacionuno-tres}) is iterative. For these cases we have $s_n=\lambda_n+\mu_n$ and $t_n=\lambda_n\mu_{n-1}$.

\medskip

Now, we give a characterization of the iterative representations when both $\{ P_n \}_{n \ge 0}$ and $\{Q_n\}_{n \ge 0}$ are SMOP's.

\begin{proposition}\label{reprminimal}
Let $\{ P_n \}_{n \ge 0}$ and $\{Q_n\}_{n \ge 0}$ be
two SMOP's  with respect to the quasi-definite linear functionals $u$ and $v$, respectively, which are related by (\ref{relacionuno-tres}). In this case, the linear functionals satisfy $(x^2+ax+b)v=(x-x_1)(x-x_2)v = k \,  u$. We denote by  $\{\beta_n\}_{n\ge0}$, $\{\gamma_n\}_{n\ge1}$ and
$\{\tilde{\beta}_n\}_{n\ge0}$, $\{\tilde{\gamma}_n\}_{n\ge1}$ the
coefficients of the three--term recurrence relation of $\{ P_n \}_{n \ge 0}$ and $\{Q_n\}_{n \ge 0}$, respectively.

Then, the representation (\ref{relacionuno-tres}) is iterative if and only if there is a sequence of complex numbers $\{ \mu_n \}_{n\ge1}$, with $\mu_1 \not= s_1$ and $\mu_n \not= 0, n\ge 1,$ such that
\begin{equation}\label{Cn}
C_n = \beta_n - \mu_{n+1} - \frac{\gamma_n}{\mu_n}=x_1 \, (\text{or} \quad x_2) \,, \quad n \ge 1 \,,
\end{equation}
and
\begin{equation}\label{must}
\mu_{n+1} = s_{n+1} - \frac{t_{n+1}}{\mu_n}, \quad n \ge 1.
\end{equation}

\end{proposition}

\textbf{Proof.} If the representation is iterative then there exists a sequence of  complex numbers $\{ \mu_n \}_{n\ge1}$, with  $\mu_n \not= 0$, such that
$R_n = P_n+ \mu_n P_{n-1}$ is an SMOP. Then, from Theorem 2 in \cite{MP},
\begin{equation*}
C_n = \beta_n - \mu_{n+1} - \frac{\gamma_n}{\mu_n}= x_1\, (\text{or} \quad x_2), \quad n \ge 1 \quad.
\end{equation*}
Substituting $R_n$ in (\ref{relacionuno-tres}), we get
\begin{equation*}
Q_n(x)=R_n(x)+(s_n-\mu_n)R_{n-1}(x)+[t_n-\mu_{n-1}(s_n-\mu_n)]P_{n-2}(x),\quad n\ge2,
\end{equation*}
and, since $Q_n$ is a linear combination of $R_n$ and $R_{n-1}$, (\ref{must}) follows.

Conversely, given $\{ \mu_n \}_{n\ge1}$ in the above conditions, from \cite[Theorem 2]{MP}, (\ref{Cn}) implies that the sequence $\{ R_n \}_{n\ge0}$ defined by $R_n = P_n+ \mu_n P_{n-1}, \quad n\ge1 ,$ is an SMOP with respect to a quasi-definite linear functional $w$ such that $(x-x_1)w = k_1 u$ where $k_1 \in \CC \setminus \{0\}$.

Taking $\lambda_n = s_n - \mu_n, \, n\ge 1$, we have  $\lambda_n \not= 0, n\ge1,$ and $t_n = \lambda_n \mu_{n-1}, n\ge2$. So, we can write
\begin{equation*}
Q_n=P_n+(\lambda_n+\mu_n)P_{n-1}+\lambda_n \mu_{n-1}P_{n-2} =
R_n + \lambda_n R_{n-1} , \quad n \ge 1 \,.
\end {equation*}
This means that (\ref{relacionuno-tres}) is iterative.
$\quad \Box$

\medskip

\noindent \textbf{Remark.} There is another constant sequence involved in this topic. Indeed, as a consequence of the above proof and according to \cite[Theorem 1]{MP}, we get
\begin{equation*}
D_n = \tilde{\beta}_n - \lambda_n - \frac{\tilde{\gamma}_{n+1}}{\lambda_{n+1}}=x_2, (\text{or} \quad x_1),  \quad n \ge 0 \,,
\end{equation*}
with $\lambda_0=0$. Moreover, we point out that, in general, the characterization of the orthogonality of a family of polynomials $\{ Q_n \}_{n\ge0}$ defined by (\ref{relacionuno-tres}) can be given (see Theorem (\ref{mainresult})) in terms of two sequences $\{ A_n \}$ and $\{ B_n \}$ which remain constant and these constant values are precisely the coefficients of the quadratic polynomial $h(x)$ which appears in the relation $h(x) v = k \,  u$. But if we have an iterative representation those sequences $\{ A_n \}$ and $\{ B_n \}$  (which in fact are very complicated) can be decomposed in two other simplest sequences $\{ C_n \}$ and $\{ D_n \}$ associated with the zeros of the polynomial $h(x)$. More precisely, when the formula (\ref{relacionuno-tres}) is iterative then $C_n + D_n = - A_n$ and $C_n D_n =B_n$.

\medskip
Observe that, from (\ref{must}), $\mu_n$ can be expressed as a finite continued fraction

\begin{equation*}
\mu_{n+1} = s_{n+1} - \cfrac{t_{n+1}}{s_n -  \cfrac{t_n}{\ddots_{\displaystyle s_2 - \cfrac{t_2}{\mu_1}}}} \quad.
\end{equation*}

\medskip
Next, as a sake of example we are going to construct a family of monic polynomials orthogonal with respect to a quasi-definite linear functional $v$  satisfying $x^2 v = x^{\alpha} e^{-x}$ and which has an iterative representation.
\medskip

\medskip

\noindent \textbf{Example. Laguerre polynomials.}

Let $\{L_n^{\alpha}\}_{n\ge0}$ be the sequence of monic Laguerre polynomials  orthogonal with respect to the positive definite linear functional $u$ defined by the weight function $x^{\alpha}\,e^{-x}$ with $ \alpha > 0$. We can take the auxiliary polynomials  $R_n=L_n^{\alpha -1}$  satisfying $R_n(x)=L_n^{\alpha}(x)+ n \,L_{n-1}^{\alpha}(x)$ (see \cite{Ch}) and the new polynomials $Q_n$ such that $Q_n(x)=R_n(x)+\lambda_n\,R_{n-1}(x)$. Then the polynomials $Q_n$ satisfy the relation (\ref{relacionuno-tres}) with $s_n=n+ \lambda_n$ and $t_n=(n-1)\,\lambda_n, \,n\ge1$. It is well known that the recurrence coefficients of $L_n^{\alpha-1}, \alpha >0,$ are $ \beta_n=2n+ \alpha, \, n \ge 0$ and $\gamma_n=n(n+ \alpha-1), \, n \ge 1$.

Using formula (\ref{Petro}) for this case, by induction we can  derive  that, for $ \alpha >0$ and $\alpha \not=1$, the values of the parameters $\lambda_n$ in terms of $\lambda_1$ are
\begin{equation}\label{mu_n}
\lambda_n=n \,\frac{\Gamma(\alpha)(\alpha-\lambda_1)+(\lambda_1-1)\frac{\Gamma(n+ \alpha)}{\Gamma(n)}}{\Gamma(\alpha)(\alpha-\lambda_1)+(\lambda_1-1)\frac{\Gamma(n-1+ \alpha)}{\Gamma(n)}}\,,
\end{equation}
and then $v$ is quasi-definite if and only if
\begin{equation*}
\Gamma(n) \Gamma(\alpha)(\alpha-\lambda_1)+(\lambda_1-1){\Gamma(n-1+ \alpha)} \not=0, \,n \ge 1.
\end{equation*}

Notice that if $\alpha \in \NN \setminus \{1\}$ then $\lambda_n$ is a rational function of $n$, namely,
$$\lambda_n=n \frac{\Gamma(\alpha)(\alpha-\lambda_1)+(\lambda_1-1)(\alpha+n-1) \dots (n+1)}{\Gamma(\alpha)(\alpha-\lambda_1)+(\lambda_1-1)(\alpha+n-2) \dots n}\,.$$

If $ \alpha=1$, then, by induction, we can also obtain,  for $ n \ge1$
\begin{equation}\label{mu_n0}
\lambda_n=n\frac{(\lambda_1-1)(1+1/2+ \dots +1/n)+1}{(\lambda_1-1)(1+1/2+ \dots +1/(n-1))+1},
\end{equation}
and  $v$ is quasi-definite  if and only if
\begin{equation*}
(\lambda_1-1)(1+1/2+ \dots +1/n)+1 \not=0, \, n\ge1.
\end{equation*}

The quasi-definite linear functional $v$ is given in terms of the Laguerre linear functional by
$$\frac{v}{v_0}=\frac{\alpha-\lambda_1}{\Gamma(\alpha)}x^{-1}\,(x^{\alpha-1}\, e^{-x})+ \delta_0$$
for $\alpha>0$.
In particular for $\alpha>1$,  we can write $$\frac{v}{v_0}=\frac{\alpha-\lambda_1}{\Gamma(\alpha)
 }\,x^{\alpha-2}\,e^{-x}+\frac{\lambda_1-1}{\alpha-1}\delta_0.$$

\bigskip

Next, as in the general case, we also present a matricial interpretation for these iterative cases.

If $w$ denotes the corresponding linear functional for $\{R_n\}_{n\ge0}$, then
$(x-x_1)w=k_1u$  where $x_1$ is a zero of $h(x)$ as well as $(x-x_2)v=k_2 w$ where $x_2$ is the other zero of $h(x)$.

It is well known (see \cite{G}) that
\begin{equation} \label {x1PnRn}
(x-x_1)P_n=R_{n+1}+r_nR_n, \quad n\ge0 \,,
\end{equation}
\noindent as well as
\begin{equation*}
(x-x_2)R_n=Q_{n+1}+r'_nQ_n, \quad n\ge0 \,,
\end{equation*}
\noindent with $r_n,r'_n \not=0$.
\medskip

As before,  ${\cal{P}}=(P_0, P_1, \dots )^T$, ${\cal{R}}=(R_0, R_1, \dots )^T$,
\noindent and  ${\cal{Q}}=(Q_0, Q_1, \dots )^T\,,$
and then the corresponding three--term recurrence relations read as
\begin{equation}\label{JPJRJQ}
x{\cal{P}}=J_{\cal{P}}\,{\cal{P}}, \quad x{\cal{R}}=J_{\cal{R}}\,{\cal{R}}, \quad x{\cal{Q}}=J_{\cal{Q}}\,{\cal{Q}}.
\end{equation}

On the other hand, from $R_n = P_n+ \mu_n P_{n-1}$ and (\ref{x1PnRn})  we have the matrix representations
\begin{equation}\label{L1U1}
{\cal{R}}=L_1{\cal{P}}, \quad (x-x_1){\cal{P}}=U_1 {\cal{R}}
\end{equation}
\noindent where $L_1$ is a lower bidiagonal matrix with $1$ as diagonal entries and $U_1$ is an upper bidiagonal matrix with $1$ as entries in the upper diagonal.

Notice that from (\ref{JPJRJQ}) and (\ref{L1U1}) we get
\begin{equation}\label{JPU1L1}
J_{\cal{P}}-x_1I=U_1L_1
\end{equation}
\noindent as well as $L_1^{-1}(J_{\cal{R}}-x_1I){\cal{R}}=U_1{\cal{R}}$. In other words
\begin{equation}\label{JRL1U1}
J_{\cal{R}}-x_1I=L_1U_1.
\end{equation}

\noindent Proceeding in a similar way
\begin{equation*}
{\cal{Q}}=L_2{\cal{R}}, \quad (x-x_2){\cal{R}}=U_2 {\cal{Q}}
\end{equation*}
\noindent where $L_2$ is a lower bidiagonal matrix with $1$ as entries in the diagonal and $U_2$ is an upper bidiagonal matrix with $1$ as entries in the upper diagonal.
Thus, we get
\begin{equation}\label{JRU2L2}
J_{\cal{R}}-x_2I=U_2L_2
\end{equation}
\noindent and
\begin{equation}\label{JQL2U2}
J_{\cal{Q}}-x_2I=L_2U_2
\end{equation}
\noindent As a consequence we can summarize the process as follows.

\noindent Step 1: Given $J_{\cal{P}}$, from $L_1$ and (\ref{JPU1L1}) we get $U_1$.

\noindent Step 2: From (\ref{JRL1U1}) we get $J_{\cal{R}}$.

\noindent Step 3: Given $J_{\cal{R}}$, from $L_2$ and (\ref{JRU2L2}) we get $U_2$.

\noindent Step 4: From (\ref{JQL2U2}) we get $J_{\cal{Q}}$.

\medskip

Notice that this is essentially the iteration of the so--called Geronimus transformation in the framework of linear spectral transformations of Stieltjes functions associated with quasi--definite linear functionals (see \cite{Z}].

Now, as in the general case,  we analyze the case of the truncated matrices.

The three--term recurrence relation for the sequence $\{ R_n \}_{n \ge 0}$ reads
\begin{equation}\label{paco4}
x\, ({\cal{R}})_n = (J_{\cal{R}})_{n+1} ({\cal{R}})_n + R_{n+1} e_{n+1} \,.
\end{equation}

Since
\begin{equation}\label{paco5}
({\cal{R}})_n = (L_1)_{n+1} ({\cal{P}})_n  \,,
\end{equation}
where $(L_1)_{n+1}$ is a lower bidiagonal matrix, namely, the $(n+1) \times (n+1)$ leading principal submatrix of $L_1$, replacing (\ref{paco5}) in (\ref{paco4})

\begin{equation*}
x \,  ({\cal{P}})_n= [(L_1)_{n+1}]^{-1}[(J_{\cal{R}})_{n+1} (L_1)_{n+1}  + \mu_{n+1} e_{n+1} e_{n+1}^T]({\cal{P}})_n + P_{n+1} e_{n+1} \,.
\end{equation*}
Thus, from (\ref{paco1}) we get $(L_1)_{n+1} (J_{\cal{P}})_{n+1} = (J_{\cal{R}})_{n+1} (L_1)_{n+1} +  \mu_{n+1} e_{n+1} e_{n+1}^T \,,$

i.e.
\begin{equation}\label{paco6}\nonumber
(J_{\cal{R}})_{n+1} =(L_1)_{n+1} [(J_{\cal{P}})_{n+1} -  \mu_{n+1} e_{n+1} e_{n+1}^T] [(L_1)_{n+1}]^{-1} \,.
\end{equation}

As a consequence, $(J_{\cal{R}})_{n+1}$ is similar to the matrix
$$(J_{\cal{P}})_{n+1} -  \mu_{n+1} e_{n+1} e_{n+1}^T$$
 that is, to the matrix obtained when we substract $\mu_{n+1}$ in the entry $(n+1) \times (n+1)$ of the Jacobi matrix $(J_{\cal{P}})_{n+1}$ and the other entries remain invariant.

On the other hand, we have
\begin{equation}\label{paco7}
(x-x_1) \, ({\cal{P}})_n = (U_1)_{n+1} ({\cal{R}})_n + R_{n+1} e_{n+1} \,,
\end{equation}
where $ (U_1)_{n+1}$ is an upper bidiagonal matrix.

According to (\ref{paco1}) and (\ref{paco5}), (\ref{paco7}) reads
\begin{equation*}
((J_{\cal{P}})_{n+1} - x_1 I_{n+1}) ({\cal{P}})_n + P_{n+1} e_{n+1}  = (U_1)_{n+1} (L_1)_{n+1} ({\cal{P}})_n + (P_{n+1} + \mu_{n+1} P_n) e_{n+1} \,.
\end{equation*}

Equivalently,
\begin{equation*}
(J_{\cal{P}})_{n+1} - \mu_{n+1} e_{n+1} e_{n+1}^T - x_1 I_{n+1} = (U_1)_{n+1} (L_1)_{n+1} \,.
\end{equation*}
But from (\ref{paco6})
\begin{equation*}
[(L_1)_{n+1}]^{-1} ((J_{\cal{R}})_{n+1} - x_1 I_{n+1}) (L_1)_{n+1} = (U_1)_{n+1} (L_1)_{n+1} \,.
\end{equation*}
Thus,
\begin{equation*}
(J_{\cal{R}})_{n+1} - x_1 I_{n+1} = (L_1)_{n+1} (U_1)_{n+1}  \,.
\end{equation*}
In a similar way,
\begin{equation*}
(J_{\cal{R}})_{n+1} - \lambda_{n+1} e_{n+1} e_{n+1}^T - x_2 I_{n+1} = (U_2)_{n+1} (L_2)_{n+1} \,,
\end{equation*}
and then
\begin{equation*}
(J_{\cal{Q}})_{n+1} - x_2 I_{n+1} = (L_2)_{n+1} (U_2)_{n+1}  \,.
\end{equation*}

\subsection{The symmetric case}

Assume that the sequence $\{ P_n \}_{n \ge 0}$ is orthogonal with respect to a symmetric linear functional $u$. Then $\beta_n=0, n\ge0$, and there exist two polynomial sequences $\{ V_n \}_{n \ge 0}$ and $\{ V_n^* \}_{n \ge 0}$ such that for all $n$,
$P_{2n}(x)=V_n(x^2)$ and $P_{2n+1}(x)=x\,V_n^*(x^2)$

It is known (see \cite{Ch}) that $\{ V_n \}_{n \ge 0}$ and $\{ V_n^* \}_{n \ge 0}$ are SMOP's with respect to the linear functionals $\sigma_u$ and $x\sigma_u$ where $(\sigma_u, x^n)=(u,x^{2n}), \, n\ge0$.

If the sequence $\{ Q_n \}_{n \ge 0}$ is defined by (\ref{relacionuno-tres}), in general, the polynomials $Q_n$ are not symmetric. More precisely, $Q_n$ are symmetric polynomials if and only if $s_n=0$ for all $n\ge1$.

Again, if $\{ Q_n \}_{n \ge 0}$ is an SMOP with respect to a symmetric linear functional $v$, then by a symmetrization process there exist $\{ R_n \}_{n \ge 0}$ and
$\{ R_n^* \}_{n \ge 0}$ SMOP's with respect to $\sigma_v$ and $x\sigma_v$, respectively, satisfying $Q_{2n}(x)=R_n(x^2)$ and $Q_{2n+1}(x)=xR_n^*(x^2)$. In this case, from (\ref{relacionuno-tres}), we have for $\ge1$
\begin{equation}\label{RnVn}
R_n(x)=V_n(x)+t_{2n}V_{n-1}(x)
\end{equation}
\begin{equation}\label{Rn*Vn*}
R_n^*(x)=V_n^*(x)+t_{2n+1}V_{n-1}^*(x)
\end{equation}
\noindent and the coefficients $t_{2n}$ and $t_{2n+1}$ can be computed using the Theorem 2 in \cite{MP}. Besides, as a consequence of the Theorem \ref{mainresult}, we get $a=0$ and the relation between the linear functionals $u$ and $v$ is: $(x^2+b)v=k \,  u$.

\medskip
Now, we construct some families of monic symmetric orthogonal polynomials which have reducible representations but not iterative representations.

\medskip

\noindent \textbf{Examples}

\medskip
\noindent \textbf{Hermite generalized polynomials.}

\medskip
Let $\{H_n^{\mu}\}_{n\ge0}$ be the sequence of monic Hermite generalized polynomials orthogonal with respect to the positive definite linear functional $u$ defined by the weight function $\vert x \vert ^{2\mu}\,e^{-x^2}$ with $ \mu > -1/2$. Consider the sequence of monic polynomials $\{Q_n\}_{n\ge0}$ orthogonal with respect to a symmetric quasi-definite functional $v$ satisfying $x^2v=k \,u$. Notice that, now the representation (\ref{relacionuno-tres}) is clearly not iterative because $xv$ is not a quasi-definite linear functional.

By a symmetrization process, $\sigma_u$ is defined by the weight function $x ^{\mu -1/2}\,e^{-x}$ and we can write $H_{2n}^{\mu}(x)=L_n^{\mu -1/2}(x^2)$, $H_{2n+1}^{\mu}(x)=x\,L_n^{\mu +1/2}(x^2)$.

From (\ref{mu_n}) and (\ref{mu_n0}) we get the explicit expressions of the parameters $t_{2n}$ in terms of $t_2$, for  $\mu >-1/2, \mu \not=1/2$, and $\mu=1/2$, respectively.

Otherwise, the relation between the two linear functionals $u$ and $v$ yields $x\sigma_v=k \,  \sigma_u$. Then the polynomials $R_n^*$ are orthogonal with respect to $x^{\mu -1/2}\,e^{-x}$ and thus $t_{2n+1}=n$ for all $n \ge1$. Moreover the quasi-definite linear functional $v$ is given in terms of the generalized Hermite linear functional $u$ by  $$\frac{v}{v_0}=\frac{\mu+1/2-t_2}{\Gamma(\mu+1/2)}x^{-2}(\vert x \vert^{2\mu}\,e^{-x^{2}})+\delta_0$$ for $\mu>-1/2$. In particular, for $\mu>1/2$ $$\frac{v}{v_0}=\frac{\mu+1/2-t_2}{\Gamma(\mu+1/2)}\,\vert x \vert^{2\mu-2}\,e^{-x^{2}}+\frac{t_2-1}{\mu-1/2}\delta_0.$$

\bigskip

\noindent \textbf{Chebyshev polynomials of the second kind.}

\bigskip

Let $\{ U_n \}_{n \ge 0}$ be the sequence of monic Chebyshev polynomials of the second kind, orthogonal with respect to the positive definite linear functional $u$ defined by the weight function $(1-x^2)^{1/2} \chi_{(-1, 1)}(x)$, with the recurrence coefficients $\beta_n=0, \,n \ge 0$ and $\gamma_n = \frac{1}{4},\, n \ge 1$. Consider the SMOP $\{Q_n\}_{n\ge0}$ orthogonal with respect to a symmetric quasi-definite linear functional $v$ satisfying $(x^2-1)v=k \,u$.

By the symmetrization process described above, we have  $\{ V_n \}_{n \ge 0}$,  $\{ R_n \}_{n \ge 0}$  are SMOP's orthogonal with respect to the linear functionals $\sigma_u$ and $\sigma_v$, respectively, where $\sigma_u$ is defined by the weight function $(1-x)^{1/2}x^{-1/2} \chi_{(0, 1)}(x)$. Besides, the recurrence coefficients for the polynomials $V_n$ are, for $n\ge1$, $\beta_n^{\sigma_u}=\gamma_{2n}+\gamma_{2n+1}=1/2$, $\gamma_n^{\sigma_u}= \gamma_{2n}\,\gamma_{2n-1}=1/16,$ and $\beta_0^{\sigma_u}=\gamma_1=1/4$.

Moreover, the relation between the two linear functionals $(x^2-1)v=k \, u$, yields  $(x-1)\sigma_v=k \,  \sigma_u$. Now, applying in (\ref{RnVn}) Theorem 2 in \cite{MP}, the parameters $t_{2n}$ satisfy $\beta_n^{\sigma_u}-t_{2n+2}-\frac{\gamma_n^{\sigma_u}}{t_{2n}}=1$ for all $n\ge1$, i.e. $t_{2n+2}=-(\frac{1}{2}+\frac{1}{16t_{2n}})$ and $k =-(\frac{3}{4}+t_2)\frac{v_0}{u_0}$. Thus, by induction, we can derive the values of parameters $t_{2n}$ in terms of $t_2.$ Indeed
\begin{equation}\label{t2n}
t_{2n}=-\frac{1}{4}\, \frac{n(4t_2+1)-1}{(n-1)(4t_2+1)-1}.
\end{equation}
On the other hand, since $(x-1)x\sigma_v=k \,  x\sigma_u$, we can apply in (\ref{Rn*Vn*}) Theorem 2 in \cite{MP}. Thus, the parameters $t_{2n+1}$ satisfy $\beta_n^{x\sigma_u}-t_{2n+3}-\frac{\gamma_n^{x\sigma_u}}{t_{2n+1}}=1$ for all $n\ge1$, i.e. $t_{2n+3}=-(\frac{1}{2}+\frac{1}{16t_{2n+1}})$ . Again, by induction, we can derive the values of parameters $t_{2n+1}$ in terms of $t_3$, and taking into account that $t_2$ and $t_3$ are related by
$t_3=\frac{12t_2+1}{4(1-4t_2)}$ we have
\begin{equation}\label{t2n+1}
t_{2n+1}=-\frac{1}{4}\, \frac{2n(4t_2+1)+(4t_2-1)}{2(n-1)(4t_2+1)+(4t_2-1)}.
\end{equation}
Moreover, if $t_2\in \CC\setminus\{0,1/2,-3/4\}$ then the linear functional $v$ is quasi-definite if and only if, for $n\ge1$
$$n(4t_2+1)-1 \not=0, \quad \text{and} \quad 2n(4t_2+1)+(4t_2-1)\not=0.$$
The expression of $v$ in terms of $u$ is  $$\frac{v}{v_0}=\frac{2}{\pi}(3/4+t_2)\,(1-x^2)^{-1/2}\chi_{(-1, 1)}(x)-(t_2+1/4)\delta_1-(t_2+1/4)\delta_{-1}.$$ See also \cite{WK2}.

\bigskip

 \noindent \textbf{Chebyshev polynomials of the first kind.}

Let $\{ T_n \}_{n \ge 0}$ be the sequence of monic Chebyshev polynomials of the first kind, orthogonal with respect to the positive definite linear functional $u$ defined by the weight function $(1-x^2)^{-1/2} \chi_{(-1, 1)}(x)$. The recurrence coefficients are $\beta_n=0, \,n \ge 0$, $\gamma_n = \frac{1}{4},\, n \ge 2$, and $\gamma_1=1/2$. Consider an SMOP $\{Q_n\}_{n\ge0}$ orthogonal with respect to a symmetric quasi-definite linear functional $v$ satisfying $(x^2-1)v=k \, u$.

Thus,  $\{ V_n \}_{n \ge 0}$,  $\{ R_n \}_{n \ge 0}$  are SMOP's  with respect to the linear functionals $\sigma_u$ and $\sigma_v$, respectively, where $\sigma_u$ is defined by the weight function $(1-x)^{-1/2}x^{-1/2} \chi_{(0, 1)}(x)$. Besides, the recurrence coefficients for the polynomials $V_n$ are, $\beta_n^{\sigma_u}=1/2$, $\gamma_n^{\sigma_u}=1/16, \, n\ge2,$ and $\gamma_1^{\sigma_u}=1/8$.

Since $(x-1)\sigma_v=k \,  \sigma_u$, applying in (\ref{RnVn}) Theorem 2 in \cite{MP}, the parameters $t_{2n}$ are
\begin{equation}
t_{2n}=-\frac{1}{4}\, \frac{n(2t_2+1)-1}{(n-1)(2t_2+1)-1}\, ,  n\ge2.
\end{equation}

Using the same arguments, from (\ref{Rn*Vn*}) we obtain
\begin{equation}
t_{2n+1}=-\frac{1}{4}\, \frac{2n(2t_2+1)+(2t_2-1)}{2(n-1)(2t_2+1)+(2t_2-1)}\, ,
\end{equation}
\noindent where we have used that $t_3=\frac{1}{4}\,\frac{6t_2+1}{1-2t_2}$.

Furthermore, if $t_2\in \CC\setminus\{0,1/2,-1/2\}$ then $v$ is quasi-definite if and only if, for $n\ge1$
$$n(2t_2+1)-1 \not=0, \quad \text{and} \quad 2n(2t_2+1)+(2t_2-1)\not=0.$$
$v$ is given in terms of $u$ by  $$\frac{v}{v_0}=\frac{(1+2t_2)}{2\pi}(1-x^2)^{-1}((1-x^2)^{-1/2}\chi_{(-1, 1)}(x))+\delta_1+\delta_{-1}.$$. See also \cite{WK1}.

\end{document}